\documentclass[12pt]{article}
\usepackage{amsmath}
\usepackage{hyperref}
\usepackage{amsfonts}
\usepackage{color}
\usepackage{amssymb}
\usepackage{graphics}
\usepackage{mathdots}
\usepackage[pdftex]{graphicx}

\def\bee{\begin{equation}}
\def\eee{\end{equation}}

\def\Li{\rm {Li}}

\def\DHLhksqrt#1#2{\setbox0=\hbox{$#1\sqrt{#2\,}$}\dimen0=\ht0
\advance\dimen0-0.2\ht0
\setbox2=\hbox{\vrule height\ht0 depth -\dimen0}%
{\box0\lower0.4pt\box2}}
\def \Li{{\rm Li}}

\allowdisplaybreaks

\begin{document}

\bigskip\bigskip\bigskip

\bigskip
\centerline{    }
\centerline{\Large\bf  Negative moments of the gaps between consecutive primes}
\bigskip\bigskip
\centerline{\large\sl Marek Wolf}
\begin{center}
Cardinal  Stefan  Wyszynski  University, Faculty  of Mathematics and Natural Sciences. College of Sciences,\\
ul. W{\'o}ycickiego 1/3,   PL-01-938 Warsaw,   Poland, e-mail:  m.wolf@uksw.edu.pl
\end{center}

\bigskip\bigskip

\bigskip\bigskip

\begin{abstract}
We derive heuristically  approximate formulas for the negative $k$--moments  $M_{-k}(x)$  of the gaps between consecutive
primes$<x $  represented  directly by  $\pi(x)$ --- the number of primes up to  $x$. In particular we propose an analytical formula
for the sum of  reciprocals of gaps between consecutive primes $<x : ~ M_{-1}(x)\sim  \frac{\pi^2(x)}{x-2\pi(x)}\log\Big(\frac{x}{2\pi(x)}\Big)
\sim  x \log \log(x)/\log^2(x)$.  We  illustrate obtained  results  by  the enormous computer data up to $x=4\times 10^{18}$.
\end{abstract}


 \bibliographystyle{abbrv}

Let $p_n$ denotes the $n$-th prime number and $d_n=p_{n+1}-p_n$ denotes the $n$-th gap between consecutive primes.
Let  us consider  the sum of  reciprocals of  gaps $d_n$  over  primes up to  $p_{n}\leq x$:
\bee
\sum_{n=2}^{\pi(x)} \frac{1}{p_{n}-p_{n-1}},
\eee
where $\pi(x)$  is,  as  usual,  the number of primes up to $x$.  To  our knowledge there is no known formula for the
above sum as a function of $x$.  We can consider   in  general  arbitrary negative moments of  $~p_{n+1}-p_n$:
\bee
M_{-k}(x) \equiv \sum_{n=2}^{\pi(x)}\frac{1}{(p_{n}-p_{n-1})^k}.
\eee
For the positive  moments
\bee
M_k(x) \equiv \sum_{n=2}^{\pi(x)}(p_{n}-p_{n-1})^k
\eee
in \cite[p.2056]{Oliveira-2013} it was  conjectured  that:
\bee
M_k(x)  \sim k! x\log^{k-1} (x).
\label{eq-Heath-Brown}
\eee
The  symbol $f(x) \sim  g(x)$ means here that $\lim_{x\rightarrow \infty} f(x)/g(x)=1$.   In  \cite {Wolf-momenty}
we predicted  the formula
\bee
M_k(x)= \frac{\Gamma(k+1)x^k}{\pi^{k-1}(x)}+\mathcal{O}_k(x).
\label{moje}
\eee
Above   $\Gamma(x)=\int_0^\infty  t^{x-1}e^{-t}dt$  is  the Gamma function generalizing the
factorial thus  the above formula  is  valid  also  for  non-integer  $k$.   By the Prime Number Theorem (PNT) the  number
of prime numbers  below $x$ is very well approximated by the logarithmic integral
\[
\pi(x)\sim  {\rm Li}(x) \equiv \int_2^x \frac{du}{\log(u)}.
\]
Integration by parts gives the asymptotic expansion which should be cut at the term $n_0=\lfloor \log(x)\rfloor$:
\bee
{\rm Li}(x) = \frac{x}{\log(x)}  + \frac{x}{\log^2(x)} + \frac{2!x}{\log^3(x)} + \frac{3!x}{\log^4(x)} +\cdots + \frac{n_0!x}{\log^{n_0+1}(x)}.
\label{PNT}
\eee
Putting in  \eqref{moje} the approximatation  $\pi(x)\sim x/\log(x)$ we recover \eqref{eq-Heath-Brown}.

Let $\tau_d(x)$ denote the number of pairs of  consecutive   primes smaller than a given bound $x$ and separated by $d$:
\bee
\tau_d(x)= \sharp\{ p_n,  p_{n+1} < x,~~{\rm  with}~ p_{n+1}-p_n=d\}.
\label{definition}
\eee
In \cite{Wolf-heuristics} (see also \cite{Wolf-1999},  \cite{Wolf-PRE}) we proposed the following formula expressing function
$\tau_d(x)$   directly by $\pi(x)$:
\bee
\tau_d(x) \sim  C_2 \prod_{p \mid d, p > 2} \frac{p - 1}{p - 2}~~ \frac{\pi^2(x)}{x}\Bigg(1-\frac{2\pi(x)}{x}\Bigg)^{{\frac{ d}{2}-1}} ~~~{\rm for} ~ d\geq 6,
\label{main}
\eee
where the  twins prime  constant
\[C_2 \equiv 2\prod_{p > 2} \biggl( 1 - \frac{1}{(p - 1)^2}\biggr) =1.320323631693739\ldots
\]
The pairs of primes separated   by $d=2$ (``twins'') and $d=4$  (``cousins'')  are special as they   always have to be
consecutive primes (with the  exception of the pair  (3,7) containing 5 in the middle)).   For $d=4$ we  adapt the expression
obtained from \eqref{main}  for $d=2$, which for $\pi(x)\sim x/\log(x)$  goes into the   the  conjecture B of G. H. Hardy
and J.E. Littlewood   \cite[eqs. (5.311) and (5.312)]{Hardy_and_Littlewood}:
\bee
\tau_2(x)\big(\approx \tau_4(x)\big)\sim  C_2 \frac{\pi^2(x)}{x} \approx C_2 \frac{x}{\log^2(x)}.
\eee
We  have
\bee
M_{-k}(x)=\sum_{n=2}^{\pi(x)}\frac{1}{(p_{n}-p_{n-1})^k} = \sum_{d=2, 4,6,\ldots} \frac{ \tau_d(x) }{d^k}.
\label{def_ujemne}
\eee
We will assume that for sufficiently regular and  decreasing  functions  $f(n)$ the  following  formula  holds:
\bee
\sum_{k=1}^\infty \prod_{p\mid k, p>2}\frac{p-1}{p-2} ~f(k) = {1 \over \prod_{p > 2}( 1 - {1 \over (p - 1)^2})} \sum_{k=1}^\infty f(k).
\label{rownosc1}
\eee
In other words we will  replace the product over $p|d$ in \eqref{main} by its mean value as  E. Bombieri and H. Davenport
\cite{Bombieri} have proved that the number  $1/\prod_{p > 2}( 1 - {1 \over (p - 1)^2}) = 2/C_2$ is the arithmetical average
of the product $\prod_{p\mid k} \frac{p-1}{p-2}$:
\bee
\frac{1}{n} \sum_{k=1}^n \prod_{p\mid k,p>2}{p-1\over p-2} =
\frac{1}{\prod_{p > 2}( 1 - \frac{1}{(p - 1)^2})} + \mathcal{O}(\log^2(n)).
\eee
Later H.L. Montgomery \cite[eq.(17.11)]{Montgomery} has improved the error term to $\mathcal{O}(\log(n))$.
Using  this trick we get further  from  \eqref{main}  and  \eqref{def_ujemne}
\bee
M_{-k}(x)\sim 2\frac{\pi^2(x)}{x-2\pi(x)}\sum_{n=1}^\infty \frac{1}{(2n)^k} \Big(1-\frac{2\pi(x)}{x}\Big)^n.
\eee
To calculate  negative  moments   we need the formula for the  series
\bee
\sum_{n=1}^\infty \frac{q^n}{n^k}\equiv {\Li}_k(q), ~~~~~|q|<1,
\eee
where  ${\Li}_k(q)$ is a  polylogarithm function of order $k$,  see, for example, \cite{Lewin-1981}  or  \cite[Sect. 25.12]{NIST-2010}.
The  ${\Li}_k(q)$   should  not be  confused  with   logarithmic integral  in \eqref{PNT},  where   ${\Li}(x)$  appears  without
any  subscript. Unfortunately  the closed  formula  for polylogarithm   is known only for   $k=1$ and is obtained by integrating
term   by term  uniformly  convergent  geometrical  series:   ${\Li}_1(q)=-\log(1-q)$. Hence we obtain
\bee
M_{-1}(x)=\sum_{n=2}^{\pi(x)} \frac{1}{p_{n}-p_{n-1}}\sim  \widetilde{M}_{-1}^{(2)}(x)\equiv \frac{\pi^2(x)}{x-2\pi(x)}\log\Big(\frac{x}{2\pi(x)}\Big).
\label{eq-minus-1}
\eee
We  use the notation $\widetilde{M}_{-k}^{(2)}(x)$  for the   analytical formula for $M_{-k}(x)$  expressed  by  $\pi(x)$,
while    $\widetilde{M}_{-k}^{(1)}(x)$   will  refer  to the  formula  for  $M_{-k}(x)$   expressed  by series in $1/\log(x)$,  see
below.   Putting here for  $\pi(x)$  a  few  first  terms from the expansion of $\Li(x)$  \eqref{PNT}  and expanding in series
of $1/\log(x)$  we obtain
\bee
M_{-1}(x) \sim  \widetilde{M}_{-1}^{(1)}(x)\equiv x\Bigg(\frac{\log(\log(x))-\log(2)}{\log^2(x)}+\frac{8\log(\log(x))-1-8\log(2)}{\log^3(x)}\Bigg).
\label{eq-minus-1a}
\eee
For  large $x$  using the notation $\log_n(x)=\log(\log_{n-1}(x))$ for the  iterated  logarithm we obtain the pleasant  formula:
\bee
\sum_{n=2}^{\pi(x)} \frac{1}{p_{n}-p_{n-1}}=1+\frac{1}{2}+\frac{1}{2}+\frac{1}{4}+\dots+\frac{1}{p_{\pi(x)}-p_{\pi(x)-1}}\sim  \frac{x\log_2(x)}{\log^2(x)}.
\label{eq-minus-1b}
\eee

During over a seven months long run of the computer  program we have collected  the values of $\tau_d(x)$ up to
$x=2^{48}\approx 2.8147\times 10^{14}$.  The data representing the function $\tau_d(x)$ were stored at values of $x$
forming the geometrical  progression with the ratio 2, i.e. at $x=2^{15}, 2^{16}, \ldots, 2^{47},
2^{48}$. Such a choice of the intermediate thresholds as powers of 2  was determined by the employed computer program
in which  the primes were coded  as bits.   The data is available for downloading from  \url{http://pracownicy.uksw.edu.pl/mwolf/gaps.zip}.
At the Tom{\'a}s   Oliveira e Silva web site \url{http://sweet.ua.pt/tos/gaps.html}  we have found values of $\tau_d(x)$ for
$x=1.61\times 10^{18}$  and  $x=4\times 10^{18}$.   In Table 1 we give a comparison of  formulas  \eqref{eq-minus-1}
and  \eqref{eq-minus-1}  with exact  values of $M_{-1}(x)$.  We  used the values of $\pi(x)$  calculated from the
identity $\sum_d \tau_d(x) = \pi(x)-1$.

\vskip 0.4cm
\begin{center}
{\sf TABLE {\bf 1}}  \\
The sum of reciprocals  of gaps between consecutive primes$<x $ compared with closed
formulas \eqref{eq-minus-1} and \eqref{eq-minus-1a}.  The  ratios initially  decrease and  next slowly  tend towards 1.
The  convergence in    second column   is  very  slow:  the ratios   are  changing  only on third places
after dot for $x$  spanning  over eleven orders.\\
\begin{tabular}{|c|c|c|} \hline
$x$ & $ M_{-1}(x)/\widetilde{M}_{-1}^{(2)}(x) $  & $  M_{-1}(x)/\widetilde{M}_{-1}^{(1)}(x)  $ \\ \hline
$2^{24}=1.6777\!  \times \! 10^{  7}$ &   0.8738&    0.7638\\ \hline
$2^{26}=6.7109\!  \times \! 10^{  7}$ &   0.8731&    0.7664\\ \hline
$2^{28}=2.6844\!  \times \! 10^{  8}$ &   0.8734&    0.7699\\ \hline
$2^{30}=1.0737\!  \times \! 10^{  9}$ &   0.8738&    0.7734\\ \hline
$2^{32}=4.2950\!  \times \! 10^{  9}$ &   0.8741&    0.7769\\ \hline
$2^{34}=1.7180\!  \times \! 10^{ 10}$ &   0.8744&    0.7803\\ \hline
$2^{36}=6.8719\!  \times \! 10^{ 10}$ &   0.8748&    0.7836\\ \hline
$2^{38}=2.7488\!  \times \! 10^{ 11}$ &   0.8751&    0.7867\\ \hline
$2^{40}=1.0995\!  \times \! 10^{ 12}$ &   0.8755&    0.7898\\ \hline
$2^{42}=4.3980\!  \times \! 10^{ 12}$ &   0.8759&    0.7927\\ \hline
$2^{44}=1.7592\!  \times \! 10^{ 13}$ &   0.8762&    0.7955\\ \hline
$2^{46}=7.0369\!  \times \! 10^{ 13}$ &   0.8766&    0.7982\\ \hline
$2^{48}=2.8147\!  \times \! 10^{ 14}$ &   0.8770&    0.8007\\ \hline
 $1.61\times 10^{18}$ &      0.8793   &    0.8145\\ \hline
 $4\times 10^{18}$  &     0.8795  &    0.8157\\ \hline

\end{tabular} \\

\end{center}
\vskip 0.4cm

For the  second negative  moment we will use the twice integrated geometrical series
\bee
\sum_{n=1}^\infty  \frac{q^n}{n(n+1)} = \frac{q+(1-q)\log(1-q)}{q}, ~~~~~~|q|<1.
\eee
as  for  large $n$  we  have $1/(n(n+1)) \approx 1/n^2$. In  this way  we obtain a crude approximation
\bee
M_{-2}(x) =\sum_{n=2}^{\pi(x)} \frac{1}{(p_{n}-p_{n-1})^2}\sim  \widetilde{M}_{-2}^{(2)}(x)\equiv \frac{1}{2}\frac{\pi^2(x)}{x-2\pi(x)}\left(1+\frac{2\pi(x)}{x-2\pi(x)}\log\Big(\frac{2\pi(x)}{x}\Big)\right)
\label{eq-minus-2}
\eee
and for $\pi(x)\sim\Li(x)$  from \eqref{PNT}    we get
\bee
M_{-2}(x)\sim  \widetilde{M}_{-2}^{(1)}(x)\equiv \frac{1}{2}\frac{ x}{\log^2(x)}\left(1-\frac{2\log(\log(x))}{\log(x)}+\mathcal{O}\Big(\frac{1}{\log^2(x)}\Big)\right).
\label{eq-minus-2a}
\eee
In Table 2 we give a comparison of these formulas with  exact  values of $M_{-2}(x)$.  The  ratios  with increasing $x$ tend
with  decreasing  speed  to 1.

\vskip 0.4cm
\begin{center}
{\sf TABLE {\bf 2}}\\
The sum of reciprocals  of squared gaps between consecutive primes$<x $ compared with closed formulas \eqref{eq-minus-2} and \eqref{eq-minus-2a}.
\begin{tabular}{|c|c|c|} \hline
$x$ & $ M_{-2}(x)/\widetilde{M}_{-2}^{(2)}(x)   $ & $  M_{-2}(x)/{\widetilde M}_{-2}^{(1)}(x)) $ \\ \hline
$2^{24}=1.6777\!  \times \! 10^{  7}$ &    1.3318 &    1.8391 \\ \hline
$2^{26}=6.7109\!  \times \! 10^{  7}$ &    1.3224 &    1.7811 \\ \hline
$2^{28}=2.6844\!  \times \! 10^{  8}$ &    1.3167 &    1.7357 \\ \hline
$2^{30}=1.0737\!  \times \! 10^{  9}$ &    1.3122 &    1.6979 \\ \hline
$2^{32}=4.2950\!  \times \! 10^{  9}$ &    1.3081 &    1.6653 \\ \hline
$2^{34}=1.7180\!  \times \! 10^{ 10}$ &    1.3044 &    1.6369 \\ \hline
$2^{36}=6.8719\!  \times \! 10^{ 10}$ &    1.3009 &    1.6119 \\ \hline
$2^{38}=2.7488\!  \times \! 10^{ 11}$ &    1.2978 &    1.5900 \\ \hline
$2^{40}=1.0995\!  \times \! 10^{ 12}$ &    1.2951 &    1.5704 \\ \hline
$2^{42}=4.3980\!  \times \! 10^{ 12}$ &    1.2926 &    1.5530 \\ \hline
$2^{44}=1.7592\!  \times \! 10^{ 13}$ &    1.2903 &    1.5373 \\ \hline
$2^{46}=7.0369\!  \times \! 10^{ 13}$ &    1.2882 &    1.5231 \\ \hline
$2^{48}=2.8147\!  \times \! 10^{ 14}$ &    1.2862 &    1.5102 \\ \hline
$1.61 \! \times  \!  10^{18}$ &    1.2767 &    1.4501 \\ \hline
$4  \! \times \! 10^{18}$ &    1.2760 &    1.4453 \\ \hline

\end{tabular} \\
\end{center}
\vskip 0.4cm

Because  the closed formula for ${\Li}_k(q)$ with $k\geq 2$  is unknown we can not obtain the conjecture for negative $k$ similar
to  \eqref{moje}.  The  approximate  formula  given in  \cite{Vepstas_2008}  is  not  convenient  for  our   purposes.
However  for  large $k$  the negative   moments are  dominated  by  smallest  gaps $d=2, 4, 6, \ldots$.  Hence  in
the defining  formula  for $M_{-k}$  we  keep  only  a  few  first  terms:
\bee
\begin{split}
M_{-k}(x)=C_2\frac{\pi^2(x)}{2^k x}\Bigg(1+\frac{1}{2^k}+\frac{2}{3^k}\left(1-\frac{2\pi(x)}{x}\right)^2+{\color{white}{\Big)} } \\
{\color{white}{\Bigg(}}+\frac{1}{4^k}\left(1-\frac{2\pi(x)}{x}\right)^3+\frac{4}{3}\frac{1}{5^k}\left(1-\frac{2\pi(x)}{x}\right)^4+\frac{2}{6^k}\left(1-\frac{2\pi(x)}{x}\right)^5+\ldots\Bigg)
\label{ujemne-k}
\end{split}
\eee
Above   we have used  explicit  values of the product $\prod_{p\mid k, p>2}\frac{p-1}{p-2}$ for $d=2, 4, 6, 8, 10$ and $12$.
For  example,  for $k=4$   keeping gaps up to $d=10$ and  developing powers of $1-2\pi(x)/x$  we obtain:
\bee
\begin{split}
M_{-4}(x)\sim { \widetilde M}_{-4}^{(2)} =C_2\frac{\pi^2(x)}{16 x}\Bigg(\frac{34081595473}{31116960000}-\frac{2500235267}{15558480000}\frac{\pi(x)}{x}+{\color{white}{\Big)}}   \\
{\color{white}{\Bigg(}} \frac{2244748963}{7779240000}\Big(\frac{\pi(x)}{x}\Big)^2-\frac{1178322017}{3889620000}\Big(\frac{\pi(x)}{x}\Big)^3+\frac{33735178}{121550625}\Big(\frac{\pi(x)}{x}\Big)^4\Bigg)
\end{split}
\eee
Putting above instead of $\pi(x)$ a  few  terms from the expansion  \eqref{PNT}  for $\Li(x)$  the   series in powers of $1/\log(x)$  follows:
\bee
\begin{split}
M_{-4}(x)\sim {\widetilde M}_{-4}^{(1)}=C_2\frac{x}{16\log^2(x)}\Big(\frac{14168273}{12960000}+\frac{14168273}{6480000}\frac{1}{\log(x)}{\color{white}{\Big)} } \\
{\color{white}{\Bigg(}}+\frac{4680091}{864000}{\Big(\frac{1}{\log(x)}\Big)^2}+ \frac{27005921}{6480000}\Big(\frac{1}{\log(x)}\Big)^3\Big)
\end{split}
\eee
In  Table 3 we  show  how  good   the above approximation  is.   We  do  not know why ratios in the third  columns are closer
to one than in second.

\bigskip

\vskip 0.4cm
\begin{center}
{\sf TABLE {\bf 3}}\\
The comparison of the formulas  for the sum of reciprocals of fourth powers of gaps between consecutive primes.
\begin{tabular}{|c|c|c|} \hline
$x$ & $ M_{-4}(x)/\widetilde{M}_{-4}^{(2)}(x) $  & $  M_{-4}(x)/\widetilde{M}_{-4}^{(1)}(x)) $ \\ \hline
$2^{24}=1.6777\!  \times \! 10^{  7}$ &  1.012003 &    1.008288 \\ \hline
$2^{26}=6.7109\!  \times \! 10^{  7}$ &  1.007536 &    1.004022 \\ \hline
$2^{28}=2.6844\!  \times \! 10^{  8}$ &  1.006260 &    1.002853 \\ \hline
$2^{30}=1.0737\!  \times \! 10^{  9}$ &  1.006038 &    1.002751 \\ \hline
$2^{32}=4.2950\!  \times \! 10^{  9}$ &  1.005621 &    1.002445 \\ \hline
$2^{34}=1.7180\!  \times \! 10^{ 10}$ &  1.005218 &    1.002190 \\ \hline
$2^{36}=6.8719\!  \times \! 10^{ 10}$ &  1.004711 &    1.001826 \\ \hline
$2^{38}=2.7488\!  \times \! 10^{ 11}$ &  1.004403 &    1.001666 \\ \hline
$2^{40}=1.0995\!  \times \! 10^{ 12}$ &  1.004104 &    1.001513 \\ \hline
$2^{42}=4.3980\!  \times \! 10^{ 12}$ &  1.003863 &    1.001417 \\ \hline
$2^{44}=1.7592\!  \times \! 10^{ 13}$ &  1.003657 &    1.001349 \\ \hline
$2^{46}=7.0369\!  \times \! 10^{ 13}$ &  1.003471 &    1.001297 \\ \hline
$2^{48}=2.8147\!  \times \! 10^{ 14}$ &  1.003311 &    1.001265 \\ \hline
$1.61 \! \times  \!  10^{18}$ &  1.002630 &    1.001258 \\ \hline
$4  \! \times \! 10^{18}$ &  1.002580 &    1.001268 \\ \hline
\end{tabular} \\
\end{center}
\vskip 0.4cm

\bigskip

We can obtain another approximate formula. Namely in \eqref{ujemne-k}  it is possible   to  sum over   all  $d$  and  separate
terms  without  $\pi(x)/x$  and with  first and second  power of $\pi(x)/x$ using the  equation  \eqref{rownosc1}.  In
this manner we obtain:
\bee
\begin{split}
M_{-k}\sim  \frac{\pi^2(x)}{2^{k-1}(x-2\pi(x))}\Bigg(\zeta(k)-\frac{2\pi(x)}{x}\Big(\zeta(k-1)-\frac{1}{2^k}\Big){\color{white}{\Bigg)} }  \\
{\color{white}{\Bigg(}}+\frac{2\pi^2(x)}{x^2}\Big(\zeta(k-2)-\zeta(k-1)-\frac{1}{2^{k-1}}~\Big) +\ldots~ \Bigg).
\end{split}
\eee
Here $\zeta(k)=\sum_n 1/n^k$  is the Riemann zeta function at integer arguments.  Above  we  have to  demand  $k\geq 4$
to  avoid infinity  of  $\zeta(1)$.

\bigskip

\noindent{\bf Acknowledgement}:  I thank Professor Roger Heath--Brown for e-mail  exchange.

\end{document}